\documentclass[11pt]{amsart}
\usepackage{amsfonts,latexsym,rawfonts,amsmath,amssymb,amsthm}
\usepackage{dsfont}
\def\be1{{\begin{equation}}}
\def\ee1{{\end{equation}}}

\def\part{\partial}

\def\ba{\begin{array}}
\def\ea{\end{array}}

\numberwithin{equation}{section}

\title{CHERN-OSSERMAN TYPE EQUALITY FOR COMPLETE SURFACES IN $\mathds{R}^{N}$}
\author{QING CHEN}
\address{University of Science and Technology of China, Department of Mathematics,
           230026 Hefei, China}
\email{qchen@ustc.edu.cn}
\author{WENJIE YANG}\address{University of Science and Technology of China, Department of Mathematics,
           230026 Hefei, China}
\email{hbxt@mail.ustc.edu.cn}
\date{}
\begin{document}

\begin{abstract}
We obtain a Chern-Osserman type equality of a complete properly immersed surface in Euclidean space, provided the $L^{2}$-norm of the second fundamental form is finite. Also, by using  a monotonicity formula, we prove that if the $L^{2}$-norm of mean curvature of a noncompact surface is finite, then it has at least quadratic area growth.
\end{abstract}

\maketitle

\section{Introduction}
Let M be a complete minimal surface in $\mathds{R}^{n}$ with finite total curvature, Chern and Osserman [2], [7] proved that
$$-\chi(M)\leq -\frac{1}{2\pi}\int_{M}K-k,\leqno{(1.1)}$$
where K is the Gauss curvature of M, $\chi(M)$ is the Euler characteristic of M and $k$ is the number of ends of M. Further results were obtained by Jorge and Meeks [5] that
$$-\chi(M)= -\frac{1}{2\pi}\int_{M}K-\lim\limits_{t\rightarrow\infty}\frac{area(M\cap B(t))}{\pi t^{2}},\leqno{(1.2)}$$
where $B(t)$ is the extrinsic ball of radius t. \\
\indent When M is a general surface properly immersed in $\mathds{R}^{n}$ with $\displaystyle\int_{M}\mid A\mid^{2} < \infty$, where A is the second fundamental form of the immersion, White [9] proved that $\displaystyle\frac{1}{2\pi}\int_{M}K$ must be an integer. In this paper, we present a general version of (1.2), where M is a general surface properly immersed in $\mathds{R}^{n}$ with the $L^{2}$-norm of the second fundamental form is finite.\\

\noindent\textbf{Theorem 1.1.} \emph{Let M be a complete properly immersed noncompact oriented surface in $\mathds{R}^{n}$, A the second fundamental form of the immersion, r the distance of  \ $\mathds{R}^{n}$ from a fixed point and $M_{t}=\{x\in M:r(x)<t\}$, $\chi(M)$ the Euler characteristic of M. Suppose $\displaystyle\int_{M}\mid A\mid^{2} < \infty$, then\\
1.$\displaystyle\lim\limits_{t \rightarrow \infty}\frac{area M_{t}}{\pi t^{2}}$ exists and is a positive integer;\\
2.$\displaystyle\lim\limits_{t \rightarrow \infty}\frac{area M_{t}}{\pi t^{2}}=\chi(M)-\frac{1}{2\pi}\int_{M}K $.}\\

 \indent Since $\displaystyle\int_{M}\mid A\mid^{2}<+\infty$, then $\displaystyle\int_{M}\mid K\mid<+\infty$ by Gauss equation. When M is a complete surface with finite total Gaussian curvature, Huber [4] proved that M has finite topological type. And Cohn-Vossen [3] obtained:
 $$2\pi\chi(M)-\int_{M}K\geq 0.\leqno{(1.3)}$$
 The explicit equality was obtained by Shiohama [8]:
$$\chi(M)-\frac{1}{2\pi}\int_{M}K=\lim\limits_{t \rightarrow \infty}\frac{D(t)}{\pi t^{2}}, \leqno{(1.4)}$$
where $D(t)$ denote the area of geodesic balls of radius t at a fixed point. Our theorem shows that (1.4) also holds with extrinsic balls instead of geodesic balls if M is properly immersed in $\mathds{R}^{n}$.\\
\indent The proof of Theorem 1.1 is based on two monotonicity formulas (Theorem 2.4). The monotonicity formulas also have an interesting application, namely, if the $L^{2}$-norm of mean curvature $H$ of the surface is finite, then it has at least quadratic area growth.\\

\noindent\textbf{Corollary 1.2.} (see  also Corollary 2.5) \emph{Let M be a complete properly immersed noncompact surface in $\mathds{R}^{n}$ with $\displaystyle\int_{M}\mid H\mid^{2}<\infty$, then the volume of the intersection of M and the extrinsic balls has at least quadratic area growth.}

\section{Preliminaries}
Let $x:M\rightarrow\mathds{R}^{n}$ be a complete properly immersed surface in $\mathds{R}^{n}$, $r$ the distance function of $\mathds{R}^{n}$ from a fixed point. For simplicity, we always assume the fixed point to be 0, unless otherwise specifie. Denote the covariant derivatve of $\mathds{R}^{n}$ and M by $\overline{{\nabla}}$ and $\nabla$ respectively. Let X, Y be two tangent vector fields of M, then
\begin{align*}\begin{split}(\overline{\nabla}^{2}r)(X,Y)&= XY(r)-\overline{{\nabla}}_{X}Y(r)\\
 &=(\nabla^{2}r)(X,Y)-\langle A(X,Y),\overline{{\nabla}}r\rangle.\end{split}\tag{2.1} \end{align*}
 The equality (2.1), together with the fact that $\displaystyle\overline{\nabla}^{2}r=\frac{1}{r}(g_{st}-dr\otimes dr)$, where $g_{st}$ denotes the standard metric of $\mathds{R}^{n}$, implies\\

 \noindent\textbf{Proposition 2.1.} \emph{For any unit tangent vector $e$ of M,
 \[(\nabla^{2}r)(e,e)=\frac{1}{r}(1-\langle e,\nabla r\rangle^{2})+\langle A(e,e),\nabla^{\perp}r\rangle,\]
 where $\nabla^{\perp}r$ is the projection of $\overline{\nabla}$ onto the normal of M.}\\

 \indent By Sard's theorem, for $a.e. t>0$, $M_{t}=\{x\in M:r(x)<t\}$ is a related compact open subset of M with the boundary $\partial M_{t}$ being a closed immersed curve of M. Let $v(t)=area M_{t}$ , A the second fundamental form of M, and $H=trA$ the mean curvature vector.\\

 \noindent\textbf{Proposition 2.2.} \emph{Suppose M is a complete properly immersed surface in $\mathds{R}^{n}$. Then for $a.e. \ t>0$,
 $$2\pi\chi(M_{t})-\int_{M_{t}}K=\frac{1}{t}\left(v'(t)+\int_{\partial M_{t}}\frac{\langle x^{\perp},H\rangle}{\mid\nabla r\mid}\right)-\int_{\partial M_{t}}\langle A(\frac{\nabla r}{\mid \nabla r\mid},\frac{\nabla r}{\mid \nabla r\mid}),\frac{\nabla^{\bot} r}{\mid \nabla r\mid} \rangle, $$
 where $x^{\perp}$ is the projection of position vector $x$ onto the normal of M.}

 \begin{proof} By the Gauss-Bonnet formula, it's sufficient to verify
 $$ \int_{\partial M_{t}}k_{g}=\frac{1}{t}(v'(t)+\int_{\partial M_{t}}\frac{\langle x^{\perp},H\rangle}{\mid\nabla r\mid})-\int_{\partial M_{t}}\langle A(\frac{\nabla r}{\mid \nabla r\mid},\frac{\nabla r}{\mid \nabla r\mid}),\frac{\nabla^{\bot} r}{\mid \nabla r\mid} \rangle, \leqno{(2.2)} $$
 where $k_{g}$ denote the geodesic curvature of $\partial M_{t}$ in M.\\
 Suppose $e$ is the unit tangent vector of $\partial M_{t}$. Since the normal of $\partial M_{t}$ is $\frac{\nabla r}{\mid\nabla r\mid}$,
 \begin{align*}\begin{split}k_{g}&=-\langle \nabla_{e}e,\frac{\nabla r}{\mid \nabla r\mid}\rangle\\
&=\frac{1}{\mid \nabla r\mid}(\nabla^{2}r)(e,e)\\
&=\frac{1}{\mid \nabla r\mid}(\frac{1}{r}+\langle A(e,e),\nabla^{\bot}r\rangle)\\
&=\frac{1}{\mid \nabla r\mid}(\frac{1}{r}+\langle H-A(\frac{\nabla r}{\mid \nabla r\mid},\frac{\nabla r}{\mid \nabla r\mid}),\nabla^{\bot}r\rangle),
\end{split}\tag{2.3}
\end{align*}
where the third equality follows by Proposition 2.1. Then by using co-area formula, $\displaystyle v'(t)=\int_{\partial M_{t}}\frac{1}{\mid\nabla r\mid}$ and the fact that $\displaystyle\nabla^{\bot}r=\frac{x^{\bot}}{r}$, we obtain (2.2).
\end{proof}

\noindent\textbf{Proposition 2.3.} \emph{Let M be a complete properly immersed surface in $\mathds{R}^{n}$, then} $$\displaystyle tv'(t)=t\int_{\partial M_{t}}\frac{\mid \nabla^{\perp} r\mid^{2}}{\mid\nabla r\mid}+2v(t)+\int_{M_{t}}\langle x^{\perp},H\rangle.$$

\begin{proof} Since $\displaystyle\frac{1}{2}\Delta r^{2}=2+\langle x,H\rangle$, integrating over $M_{t}$ and using the Green's formula,
$$t\int_{\partial M_{t}}\mid\nabla r\mid=2v(t)+\int_{M_{t}}\langle x,H\rangle.\leqno{(2.4)}$$
By the co-area formula ,
$$ v'(t)=\int_{\partial M_{t}}\frac{1}{\mid\nabla r\mid}.$$
So we have,
\begin{align*} tv'(t)&=t\left(\int_{\partial M_{t}}\frac{1}{\mid\nabla r\mid}-\int_{\partial M_{t}}\mid\nabla r\mid\right)+t\int_{\partial M_{t}}\mid\nabla r\mid\\
 &=t\int_{\partial M_{t}}\frac{\mid \nabla^{\perp} r\mid^{2}}{\mid\nabla r\mid}+2v(t)+\int_{M_{t}}\langle x^{\perp},H\rangle.\qedhere \end{align*}
 \end{proof}

 \noindent\textbf{Theorem 2.4.} \emph{Let M be a complete properly immersed surface in $\mathds{R}^{n}$, r the distance of $\mathds{R}^{n}$ from a fixed point $x_{0}$, H the mean curvature of M, $M_{t}=\{x\in M:r(x)<t\}$, $v(t)=area{M_{t}}$, then both
 $$u_{1}(t)\triangleq\frac{v(t)}{t^{2}}-\frac{1}{2t^{2}}\int_{M_{t}}\mid (x-x_{0})^{\perp}\mid\mid H\mid+\frac{1}{16}\int_{M_{t}}\mid H\mid^{2}$$
 and
 $$u_{2}(t)\triangleq\frac{v(t)}{t^{2}}-\frac{1}{t^{2}}\int_{M_{t}}\mid (x-x_{0})^{\perp}\mid\mid H\mid+\frac{1}{4}\int_{M_{t}}\mid H\mid^{2}$$
 are monotone nondecreasing in t.}

 \begin{proof} For simplicity, we assume $x_{0}=0$. By Proposition 2.3, we have
$$tv'(t)\geq t\int_{\partial M_{t}}\frac{\mid \nabla^{\perp} r\mid^{2}}{\mid\nabla r\mid}+2v(t)-\int_{M_{t}}\mid x^{\perp}\mid\mid H\mid.  \leqno{(2.5)}$$
By co-area formula and the weighted mean value inequalities,
\begin{align}\begin{split}\frac{d}{dt}(\int_{M_{t}}\mid x^{\perp}\mid\mid H\mid)&=\int_{\partial M_{t}}\frac{\mid x^{\perp}\mid\mid H\mid}{\mid\nabla r\mid}\\
&= t\int_{\partial M_{t}}\frac{\mid \nabla^{\perp}r\mid\mid H\mid}{\mid\nabla r\mid}\\
&\leq 2\int_{\partial M_{t}}\frac{\mid \nabla^{\perp} r\mid^{2}}{\mid\nabla r\mid}+\frac{t^{2}}{8}\int_{\partial M_{t}}\frac{\mid H\mid^{2}}{\mid\nabla r\mid}.
\end{split}\tag{2.6}\end{align}
Combining (2.5) and (2.6), we have
$$tv'(t)\geq \frac{t}{2}\left((\int_{M_{t}}\mid x^{\perp}\mid\mid H\mid)'-\frac{t^{2}}{8}\int_{\partial M_{t}}\frac{\mid H\mid^{2}}{\mid\nabla r\mid}\right)+2v(t)-\int_{M_{t}}\mid x^{\perp}\mid\mid H\mid,\leqno{(2.7)}$$
or equivalently,
$$tv'(t)-2v(t)-\frac{1}{2}\left(t(\int_{M_{t}}\mid x^{\perp}\mid\mid H\mid)'-2\int_{M_{t}}\mid x^{\perp}\mid\mid H\mid\right)+\frac{t^{3}}{16}\int_{\partial M_{t}}\frac{\mid H\mid^{2}}{\mid\nabla r\mid}\geq 0.\leqno{(2.8)}$$
Dividing both sides of (2.8) by $t^{3}$ yields
$$\frac{d}{dt}\left(\frac{v(t)}{t^{2}}-\frac{1}{2}\frac{\int_{M_{t}}\mid x^{\perp}\mid\mid H\mid}{t^{2}}+\frac{1}{16}\int_{M_{t}}\mid H\mid^{2}\right)\geq 0,\leqno{(2.9)}$$
this proves that $u_{1}(t)$ is monotone nondecreasing in t.\\
\indent If we make slight modifications to (2.5) and (2.6), we have
$$tv'(t)\geq t\int_{\partial M_{t}}\frac{\mid \nabla^{\perp} r\mid^{2}}{\mid\nabla r\mid}+2v(t)-2\int_{M_{t}}\mid x^{\perp}\mid\mid H\mid,  \leqno{(2.5)'}$$
and
$$\frac{d}{dt}(\int_{M_{t}}\mid x^{\perp}\mid\mid H\mid)\leq \int_{\partial M_{t}}\frac{\mid \nabla^{\perp} r\mid^{2}}{\mid\nabla r\mid}+\frac{t^{2}}{4}\int_{\partial M_{t}}\frac{\mid H\mid^{2}}{\mid\nabla r\mid}. \leqno{(2.6)'}$$
Combining $(2.5)'$ and $(2.6)'$, we obtain
$$\frac{d}{dt}\left(\frac{v(t)}{t^{2}}-\frac{\int_{M_{t}}\mid x^{\perp}\mid\mid H\mid}{t^{2}}+\frac{1}{4}\int_{M_{t}}\mid H\mid^{2}\right)\geq 0,\leqno{(2.9)'} $$
i.e.  $u_{2}(t)$ is monotone nondecreasing in t.
\end{proof}

\noindent\emph{Remark 2.4.} From the poof, we can see that the theorem is also valid for noncomplete surface, for t with $\partial M\cap B_{x_{0}}(t)=\emptyset$, where $B_{x_{0}}(t)$ is the ball in $\mathds{R}^{n}$ of radius t and centered at $x_{0}$. \\

\indent By Theorem 2.4, we can get various volume estimates under suitable restrictions on mean curvature $H$.\\

\noindent\textbf{Corollary 2.5.} \emph{Let M be a complete properly immersed noncompact surface in $\mathds{R}^{n}$ with $\displaystyle\int_{M}\mid H\mid^{2}<\infty$, then the volume of the intersection of M and the extrinsic balls has at least quadratic area growth.}

\begin{proof} Without loss of generality, we assume the center of the extrinsic balls to be 0. Since $\displaystyle\int_{M}\mid H\mid^{2}<\infty$, for a given $\varepsilon>0$, there exists $R>0$, such that
$$\int_{M\setminus B_{0}(R)}\mid H\mid^{2}<\varepsilon.$$
Now for $t>R$ large enough, choosing a point $p\in M\cap \partial B_{0}(\frac{t+R}{2})$, then $B_{p}(\frac{t-R}{2})\subset B_{0}(t)\setminus B_{0}(R)$, so we have
$$\int_{M\cap B_{p}(\frac{t-R}{2})}\mid H\mid^{2}<\varepsilon,\ \ Vol(M\cap B_{0}(t))\geq Vol(M\cap B_{p}(\frac{t-R}{2})).\leqno{(2.10)}$$
Taking $x_{0}=p$ in Theorem 2.4, then we have
$$u_{1}(\frac{t-R}{2})\geq \lim\limits_{r\rightarrow 0}u_{1}(r)=\pi.\leqno{(2.11)}$$
Combining (2.10) and (2.11), we obtain
\begin{align*}Vol(M\cap B_{0}(t))&\geq Vol(M\cap B_{p}(\frac{t-R}{2}))\\
&\geq \frac{(t-R)^{2}}{4}\left(u_{1}(\frac{t-R}{2})-\frac{1}{16}\int_{M_{\frac{t-R}{2}}}\mid H\mid^{2}\right)\\
&\geq \frac{\pi-\varepsilon}{4}(t-R)^{2}.\end{align*}
The conclusion follows by chosing $\varepsilon$ small.
\end{proof}

\section{Proof of Theorem 1.1}

\noindent\textbf{Lemma 3.1.} \emph{Let M be as in Theorem 1.1, then both $\displaystyle\lim\limits_{t \rightarrow \infty}\frac{v(t)}{t^{2}}$ and $\displaystyle\lim\limits_{t \rightarrow \infty}\frac{\int_{M_{t}}\mid x^{\perp}\mid\mid H\mid}{t^{2}}$ exist.}

\begin{proof} First we prove:\\
\noindent\textbf{Claim:} $\displaystyle\liminf \limits_{t\rightarrow \infty}\frac{v}{t^{2}}<+\infty.$ \\
\textbf{Proof of the claim:} Since by the weighted mean value inequality, $$\mid\frac{1}{t}\int_{\partial M_{t}}\frac{\langle x^{\perp},H\rangle}{\mid\nabla r\mid}\mid\leq \int_{\partial M_{t}}\frac{\mid H\mid}{\mid\nabla r\mid}\leq \frac{1}{2}\left(t\int_{\partial M_{t}}\frac{\mid H\mid^{2}}{\mid\nabla r\mid}+\frac{v'}{t}\right),\leqno{(3.1)}$$ by Proposition 2.2, we have
\begin{align}\begin{split} 2\pi \chi(M_{t})-\int\limits_{M_{t}}K&\geq \frac{1}{t}\left(v'(t)-\mid\int_{\partial M_{t}}\frac{\langle x^{\perp},H\rangle}{\mid\nabla r\mid}\mid\right)-\int\limits_{\partial M_{t}}\mid A\mid \frac{\mid \nabla^{\bot} r\mid }{\mid \nabla r\mid}\\
&\geq \frac{v'}{2t}-\frac{t}{2}\int_{\partial M_{t}}\frac{\mid H\mid^{2}}{\mid\nabla r\mid}-\int\limits_{\partial M_{t}}(\frac{t\mid A\mid^{2}}{2\mid \nabla r\mid}+\frac{\mid \nabla^{\bot} r\mid^{2} }{2t\mid \nabla r\mid})\\
(by \ Proposition \ 2.3)\ \ &= \frac{v'}{2t}-\frac{t}{2}\int_{\partial M_{t}}\frac{\mid H\mid^{2}}{\mid\nabla r\mid}-\frac{t}{2}\int_{\partial M_{t}}\frac{\mid A\mid^{2}}{\mid\nabla r\mid}-\frac{tv'-(2v(t)+\int_{M_{t}}\langle x^{\perp},H\rangle)}{2t^{2}}\\
&=\frac{v}{t^{2}}-\frac{t}{2}\int_{\partial M_{t}}\frac{\mid H\mid^{2}}{\mid\nabla r\mid}-\frac{t}{2}\int_{\partial M_{t}}\frac{\mid A\mid^{2}}{\mid\nabla r\mid}+\frac{\int_{M_{t}}\langle x^{\perp},H\rangle}{2t^{2}}\\
&\geq \frac{v}{t^{2}}-\frac{t}{2}\int_{\partial M_{t}}\frac{\mid H\mid^{2}}{\mid\nabla r\mid}-\frac{t}{2}\int_{\partial M_{t}}\frac{\mid A\mid^{2}}{\mid\nabla r\mid}-\frac{\sqrt{v\int_{M_{t}}\mid H\mid^{2}}}{2t}\\
&\geq \frac{v}{2t^{2}}-\frac{t}{2}\int_{\partial M_{t}}\frac{\mid H\mid^{2}}{\mid\nabla r\mid}-\frac{t}{2}\int_{\partial M_{t}}\frac{\mid A\mid^{2}}{\mid\nabla r\mid}-\frac{1}{8}\int_{M_{t}}\mid H\mid^{2},\end{split}\tag{3.2}
\end{align}
where we use the weighted mean value inequalities in the second and the last equality, while the second equality count backwards follows from Cauchy's inequality.\\
Since $\displaystyle\int_{M}\mid H\mid^{2}+\int_{M}\mid A\mid^{2}<+\infty$, there exists a sequence $\{\tau_{i}\}$ diverging to infinity such that
$$t(\int_{\partial M_{t}}\frac{\mid H\mid^{2}}{\mid\nabla r\mid}+\int_{\partial M_{t}}\frac{\mid A\mid^{2}}{\mid\nabla r\mid})\bigg|_{t=\tau_{i}}\longrightarrow 0\ \ \  as \ \ \ i\longrightarrow \infty.\leqno{(3.3)}$$ Otherwise, we  must have $\displaystyle\liminf\limits_{t\rightarrow\infty}t(\int_{\partial M_{t}}\frac{\mid H\mid^{2}}{\mid\nabla r\mid}+\int_{\partial M_{t}}\frac{\mid A\mid^{2}}{\mid\nabla r\mid})=\delta>0$. \ So for sufficient large t, we have$$\displaystyle t(\int_{\partial M_{t}}\frac{\mid H\mid^{2}}{\mid\nabla r\mid}+\int_{\partial M_{t}}\frac{\mid A\mid^{2}}{\mid\nabla r\mid})>\frac{\delta}{2},$$  i.e. $$\displaystyle\int_{\partial M_{t}}\frac{\mid H\mid^{2}}{\mid\nabla r\mid}+\int_{\partial M_{t}}\frac{\mid A\mid^{2}}{\mid\nabla r\mid}>\frac{\delta}{2t}.$$When you integrate t, by the co-area formula, it is as bounded on the left as it is diverging on the right, a contradiction.\\
Then taking $t=\tau_{i}$ in (3.2), together with the fact that \\
\[\displaystyle\chi(M_{t})\leq 1,\ \ \displaystyle\left|\int_{M_{t}}K\right|\leq\frac{1}{2}\int_{M}\mid A\mid^{2}<+\infty  \ \  and  \ \ \displaystyle\int_{M_{t}}\mid H\mid^{2}\leq 2\int_{M}\mid A\mid^{2}<+\infty,\]
  we have $\displaystyle\limsup\limits_{i\rightarrow \infty}\frac{v(t_{i})}{t_{i}^{2}}<+\infty$, which implies
$\displaystyle\liminf \limits_{t\rightarrow \infty}\frac{v}{t^{2}}\leq\limsup\limits_{i\rightarrow \infty}\frac{v(t_{i})}{t_{i}^{2}}<+\infty$. This proves the claim.\\
\indent Let $u_{1}(t)$ and $u_{2}(t)$ be as in Theorem 2.4 with $x_{0}=0$. By the claim, we have\begin{align*}\begin{split}\liminf\limits_{t\rightarrow\infty}u_{1}(t)&\leq \liminf\limits_{t\rightarrow\infty}\frac{v(t)}{t^{2}}+\frac{1}{16}\int_{M}\mid H\mid^{2}<+\infty,\\
\liminf\limits_{t\rightarrow\infty}u_{2}(t)&\leq \liminf\limits_{t\rightarrow\infty}\frac{v(t)}{t^{2}}+\frac{1}{4}\int_{M}\mid H\mid^{2}<+\infty.
\end{split} \tag{3.4}\end{align*}
Combining (3.4) and Theorem 2.4, we know that both $u_{1}(t)$ and $u_{2}(t)$ have finite limit as $t\rightarrow\infty$.\\
Since\begin{align*}\begin{split}\frac{v(t)}{t^{2}}&=2u_{1}(t)-u_{2}(t)+\frac{1}{8}\int_{M_{t}}\mid H\mid^{2},\\
\frac{\int_{M_{t}}\mid x^{\perp}\mid\mid H\mid}{t^{2}}&=2u_{1}(t)-2u_{2}(t)+\frac{3}{8}\int_{M_{t}}\mid H\mid^{2},\end{split} \tag{3.5}\end{align*}
we conclude that both $\displaystyle\lim\limits_{t \rightarrow \infty}\frac{v(t)}{t^{2}}$ and $\displaystyle\lim\limits_{t \rightarrow \infty}\frac{\int_{M_{t}}\mid x^{\perp}\mid\mid H\mid}{t^{2}}$ exist.
\end{proof}

\noindent\textbf{Lemma 3.2.} \emph{There exists a sequence $\{t_{k}\}$ diverging to infinity such that}
\begin{align*}&\lim\limits_{k \rightarrow \infty}\frac{v'(t_{k})}{t_{k}}=\lim\limits_{k \rightarrow \infty}\frac{2v(t_{k})}{t_{k}^{2}}, \lim\limits_{k \rightarrow \infty}\frac{2}{t_{k}^{2}}\int_{M_{t_{k}}}\mid x^{\perp}\mid\mid H\mid=\lim\limits_{k \rightarrow \infty}\frac{1}{t_{k}}\int_{\partial M_{t_{k}}}\frac{\mid x^{\perp}\mid\mid H\mid}{\mid\nabla r\mid}=0,\tag{i}\\
 &\lim\limits_{k \rightarrow \infty}t_{k}\int_{\partial M_{t_{k}}}\frac{\mid H\mid^{2}}{\mid\nabla r\mid}=0,\ \lim\limits_{k \rightarrow \infty}t_{k}\int_{\partial M_{t_{k}}}\frac{\mid A\mid^{2}}{\mid\nabla r\mid}=0.\tag{ii}\end{align*}

\begin{proof} Let $u_{1}(t)$, $u_{2}(t)$ be as in Lemma 3.1. Since $u_{1}(t)+u_{2}(t)+\int_{M_{t}}\mid H\mid^{2}+\int_{M_{t}}\mid A\mid^{2}$ is bounded, arguing as in the proof of the claim in Lemma 3.1, we know that there is a sequence $\{t_{k}\}$ diverging to infinity such that
$$t\frac{d}{dt}\left(u_{1}(t)+u_{2}(t)+\int_{M_{t}}\mid H\mid^{2}+\int_{M_{t}}\mid A\mid^{2}\right)\bigg|_{t=t_{k}}\longrightarrow 0\ \  as\ \  k\longrightarrow \infty.\leqno{(3.6)}$$
Since derivative of each function in left side of (3.6) is nonnegative, we have
 \begin{align*}\begin{split}&t_{k}u'_{1}(t_{k})\rightarrow 0,t_{k}u'_{2}(t_{k})\rightarrow 0 \ \ and\\&t\left(\int_{M_{t}}\mid H\mid^{2}\right)'\bigg|_{t=t_{k}}\rightarrow 0,t\left(\int_{M_{t}}\mid A\mid^{2}\right)'\bigg|_{t=t_{k}}\rightarrow 0\end{split}\tag{3.7}\end{align*}as\ \  $k\rightarrow \infty.$
 Combining (3.5) and (3.7), we get
 \begin{align}\begin{split}&t\left(\frac{v(t)}{t^{2}}\right)'\bigg|_{t=t_{k}}\rightarrow 0,t\left(\frac{\int_{M_{t}}\mid x^{\perp}\mid\mid H\mid}{t^{2}}\right)'\bigg|_{t=t_{k}}\rightarrow 0 \ \ and\\&t\left(\int_{M_{t}}\mid H\mid^{2}\right)'\bigg|_{t=t_{k}}\rightarrow 0,t\left(\int_{M_{t}}\mid A\mid^{2}\right)'\bigg|_{t=t_{k}}\rightarrow 0\end{split}\tag{3.8}\end{align}
 as\ \  $k\rightarrow \infty.$
 So we obtain
 \begin{align}\begin{split}&\lim\limits_{k \rightarrow \infty}\frac{v'(t_{k})}{t_{k}}=\lim\limits_{k \rightarrow \infty}\frac{2v(t_{k})}{t_{k}^{2}}, \
 \lim\limits_{k \rightarrow \infty}\frac{1}{t_{k}}\left(\int_{M_{t_{k}}}\mid x^{\perp}\mid\mid H\mid\right)'=\lim\limits_{k \rightarrow \infty}\frac{2}{t_{k}^{2}}\int_{M_{t_{k}}}\mid x^{\perp}\mid\mid H\mid \ \ and\\
 &\lim\limits_{k \rightarrow \infty}t_{k}\int_{\partial M_{t_{k}}}\frac{\mid H\mid^{2}}{\mid\nabla r\mid}=0,\ \lim\limits_{k \rightarrow \infty}t_{k}\int_{\partial M_{t_{k}}}\frac{\mid A\mid^{2}}{\mid\nabla r\mid}=0,\end{split}\tag{3.9}\end{align}
 where we use the fact that $\lim\limits_{t \rightarrow \infty}\frac{v(t)}{t^{2}}$ and $\lim\limits_{t \rightarrow \infty}\frac{\int_{M_{t}}\mid x^{\perp}\mid\mid H\mid}{t^{2}}$ exist by Lemma 3.1, this proves (ii).\\
 By co-area formula, when $k\longrightarrow\infty$,
 \begin{align}\begin{split}\frac{1}{t}\frac{d}{dt}\left(\int_{M_{t}}\mid x^{\perp}\mid\mid H\mid\right)\bigg|_{t=t_{k}}&=\frac{1}{t_{k}}\int_{\partial M_{t_{k}}}\frac{\mid x^{\perp}\mid\mid H\mid}{\mid\nabla r\mid}\\
 &\leq \int_{\partial M_{t_{k}}}\frac{\mid H\mid}{\mid\nabla r\mid}\\
 &\leq \sqrt{v'(t_{k})\int_{\partial M_{t_{k}}}\frac{\mid H\mid^{2}}{\mid\nabla r\mid}}\\
 &=\sqrt{\frac{v'(t_{k})}{t_{k}}t_{k}\int_{\partial M_{t_{k}}}\frac{\mid H\mid^{2}}{\mid\nabla r\mid}}\\
 &\longrightarrow 0.\end{split}\tag{3.10}
 \end{align}
 Combining (3.9) and (3.10), we have
 $$\lim\limits_{k \rightarrow \infty}\frac{2}{t_{k}^{2}}\int_{M_{t_{k}}}\mid x^{\perp}\mid\mid H\mid=\lim\limits_{k \rightarrow \infty}\frac{1}{t_{k}}\int_{\partial M_{t_{k}}}\frac{\mid x^{\perp}\mid\mid H\mid}{\mid\nabla r\mid}=0.\leqno({3.11)}$$
 Then (i) follows from (3.9) and (3.11).
 \end{proof}

 \noindent\textbf{Proof of  Theorem 1.1}\ \  By Proposition 2.3, we have
 \begin{align*}\begin{split}\bigg|\int_{\partial M_{t}}\langle A(\frac{\nabla r}{\mid \nabla r\mid},\frac{\nabla r}{\mid \nabla r\mid}),\frac{\nabla^{\bot} r}{\mid \nabla r\mid} \rangle\bigg|&\leq \int\limits_{\partial M_{t}}\mid A\mid \frac{\mid \nabla^{\bot} r\mid }{\mid \nabla r\mid}\\
 &\leq \int\limits_{\partial M_{t}}(\frac{t\mid A\mid^{2}}{2\mid \nabla r\mid}+\frac{\mid \nabla^{\bot} r\mid^{2} }{2t\mid \nabla r\mid})\\
 &=\frac{t}{2}\int_{\partial M_{t}}\frac{\mid A\mid^{2}}{\mid\nabla r\mid}+\frac{tv'-(2v(t)+\int_{M_{t}}\langle x^{\perp},H\rangle)}{2t^{2}},\end{split} \tag{3.12}
 \end{align*}
 then Lemma 3.2 implies
 $$\lim\limits_{k\rightarrow\infty}\bigg|\int_{\partial M_{t_{k}}}\langle A(\frac{\nabla r}{\mid \nabla r\mid},\frac{\nabla r}{\mid \nabla r\mid}),\frac{\nabla^{\bot} r}{\mid \nabla r\mid} \rangle\bigg|=0.\leqno{(3.13)}$$
 Taking $t=t_{k}$ in Proposition 2.2 and letting $k\rightarrow\infty$, together with (3.13) and Lemma 3.2, we get
 $$2\pi\lim\limits_{k\rightarrow\infty}\chi(M_{t_{k}})-\int_{M}K=\lim\limits_{k\rightarrow\infty}\frac{2v(t_{k})}{t_{k}^{2}},\leqno{(3.14)}$$
 which implies
 $$\lim\limits_{t\rightarrow\infty}\frac{2v(t)}{t^{2}}\leq 2\pi\chi(M)-\int_{M}K.\leqno{(3.15)}$$
 Since the extrinsic distance is smaller than intrinsic distance, we clearly have
 $$\lim\limits_{t\rightarrow\infty}\frac{v(t)}{t^{2}}\geq \lim\limits_{t\rightarrow\infty}\frac{D(t)}{t^{2}},\leqno{(3.16)}$$
 where D(t) is the area of geodesic balls of radius t at a fixed point.\\
 Combining (1.4), (3.15) and (3.16), we conclude that
 $$\lim\limits_{t\rightarrow\infty}\frac{2v(t)}{t^{2}} = 2\pi\chi(M)-\int_{M}K.\leqno{(3.17)}$$
 Furthermore, by the main theorem of White [9], we know that $\frac{1}{2\pi}\int_{M}K$ is an integer, so is $\lim\limits_{t\rightarrow\infty}\frac{v(t)}{\pi t^{2}}$, and this limit must be positive by Corollary 2.5. This completes the proof of Theorem 1.1.\\

\noindent\textbf{Corollary 3.3}\emph{ Let M be a complete properly immersed noncompact oriented surface in $\mathds{R}^{n}$ with $\displaystyle\int_{M}\mid A\mid^{2} <4\pi$, then $\chi(M)=1$.}

 \renewcommand\refname{R\scriptsize EFERENCES}

\end{document}